\documentclass[12pt]{article}
\usepackage{amsmath} 
\usepackage{graphicx}
\usepackage[all]{xy}
\usepackage{latexcad}
\usepackage{amssymb}
\input diagxy
\usepackage[usenames]{color}
\newtheorem{definition}{Definition}[section]

\newtheorem{prop}{Proposition}[section]

\newtheorem{rem}{Remark}[section]

\newtheorem{thm}{Theorem}[section]

\newtheorem{dfn}{Definition}[section]

\def\mathopdef#1{\expandafter\def\csname#1\endcsname{\mathop{\rm#1}\nolimits}}
\def\mathoplsdef#1{\expandafter\def\csname#1\endcsname{\mathop{\rm#1}}}
\def\mathbfdef#1{\expandafter\def\csname#1\endcsname{{\rm\bf#1}}}
\def\mathrmdef#1{\expandafter\def\csname#1\endcsname{{\rm#1}}}

\newcommand{\tn}{\otimes}

\newcommand{\pf}{\vspace{\baselineskip}\noindent{\bf Proof. }}
\newcommand{\qed}{\hfill\rule{4pt}{8pt}\par\vspace{\baselineskip}}

\newcommand{\bB}{{\bf B}}
\newcommand{\bC}{{\bf C}}

\newcommand{\BB}{\mathbf{BB}}
\newcommand{\Z}{\mathbb{Z}}

\newcommand{\eps}{{\epsilon}}

\newcommand{\tw}{{\tau}}

\newcommand{\s}{{\sigma}}

\newcommand{\rel}{\mathbf{Rel}}

\newcommand{\tr}{\mathbf{TRel}}

\newcommand{\tang}{\mathbf{Tangle}}

\title{Blocked-braid Groups}
\author{D. Maglia\\
Dipartimento di Scienza e Alta Tecnologia\\
Universit\`a dell' Insubria\\
via Carloni, 78, Como, Italy\\
{\tt dvd.maglia@libero.it}
\and
N. Sabadini\\
Dipartimento di Scienza e Alta Tecnologia\\
Universit\`a dell' Insubria\\
via Carloni, 78, Como, Italy\\
{\tt nicoletta.sabadini@uninsubria.it}
\and
R. F. C. Walters\\
Dipartimento di Scienza e Alta Tecnologia\\
Universit\`a dell' Insubria\\
via Carloni, 78, Como, Italy\\
{\tt robert.walters@uninsubria.it}}

\begin{document}
\maketitle
\begin{abstract}
   We introduce and study a family of groups $\BB_n$, called the blocked-braid groups,
 which are quotients of 
Artin's braid groups $\bB_n$, and have the corresponding symmetric groups $\Sigma_n$ as quotients. They are defined by adding a certain class of geometrical modifications to braids.  They arise in the study of commutative Frobenius algebras and tangle algebras in braided strict monoidal categories.  A fundamental equation true in $\BB_n$ is Dirac's Belt Trick; that torsion through $4\pi$ is equal to the identity. We show that $\BB_n$ is finite for $n=1,2$ and $3$ but infinite for $n>3$. 
\end{abstract}

\newpage

\section{Introduction}\label{sec-intro}
The notion of \emph{tangled circuit diagram} was introduced in \cite{RSW-TangCirc-2012} as part of a programme (\cite{KSW97b,KSW00a, KSW00b, KSW02,RSW04,RSW05,RSW08,dFSWcsg}) to study the kinds of parallel and sequential networks which occur in Computer Science. Allowing tangling introduced questions of geometric interest, and the search for invariants lead to connections with knot invariants.

This paper is an attempt to concentrate attention on some very specific tangled circuits which we called blocked braids, which consist of two components with the same number of ports joined by a braid of wires. Such circuits have a group structure (analogous to the sum of knots) and we describe some progress in understanding these groups. 

This work appeared in the Laurea Magistrale Thesis of Davide Maglia \cite{M12}.

\section{Blocked-braid groups}

\begin{dfn}\label{def-braided-monoidal}
 A \emph{braided strict monoidal category} (\cite{jsbmc}) is a category $\bC$ with a functor, called tensor,
$\tn : \bC \times \bC \to \bC$ and a ``unit'' object $I$
together with a natural family of isomorphisms $\tw_{A,B} : A \tn B \to B \tn A$ called twist
satisfying
\begin{itemize}
\item[1)] $\tn$ is associative and unitary on objects and arrows,
\item[2)] the following diagrams commute for objects $A, B, C$: 
$$
\bfig
\place(-400,200)[B1:]
\Vtriangle/>`>`<-/<600,400>[A\tn B\tn C`B\tn C\tn A`B\tn A\tn C;\tw`\tw\tn 1`1\tn\tw]
\efig
$$
and
$$
\bfig
\place(-400,200)[B2:]
\Vtriangle/>`>`<-/<600,400>[A\tn B\tn C`C\tn A\tn B`A\tn C\tn B;\tw`1\tn\tw`\tw\tn 1]
\efig
$$
\end{itemize}
We will denote the n-fold tensor product of an object $A$ by $A^n$.
\end{dfn}

\begin{dfn}\label{def-tangle-algebra}
 A tangle algebra in a braided strict monoidal category (with twist $\tau$) is an object $A$ equipped with 
arrows $\eta : I \to A \tn A$ and $\eps : A \tn A \to I$ that satisfy the equations:
\begin{itemize}
\item[(i)] $(\eps\tn 1_A) (1_A\tn \eta) = 1_A = (1_A\tn \eps)(\eta\tn 1_A) $
\item[(ii)]  $\eps \tw_{A,A} = \eps$ and $\tw_{A,A}\eta  = \eta$.
\end{itemize}
\end{dfn}

\begin{dfn}\label{def-tangle-category}\cite{FY89}
The category $\tang_n$ is the free braided strict monoidal category generated by an object $A$ with a tangle algebra structure, and two arrows $R:I\to A^n$ and $S:A^n\to I$. 
\end{dfn}

\begin{dfn}\label{def-braid}
We call an arrow $B:A^n\to A^n$ in $\tang_n$ a \emph{braid in $\tang_n$} if it is a composite of arrows each of which is a tensor of identities and arrows $\tw_{A,A}$.
It is straightforward to see that since $\tw$ has an inverse then so also has any braid in $\tang_n$.
\end{dfn}
We can picture braids in an obvious way (though notice that the order of composition of arrows in the expression is reversed in the picture). For example, the braid 
$$ (1\tn \tw)(\tw\tn 1)(1\tn \tw)(\tw\tn 1)(1\tn \tw)(\tw\tn 1)$$
 is pictured as:
 \begin{center}
 \includegraphics[width=250mm,bb=-100 0 900 100]{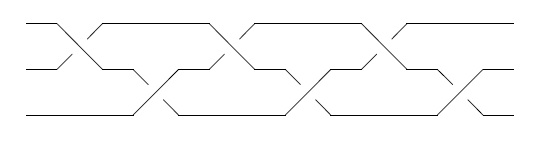}
 \end{center}

\begin{dfn}\label{def-blocked-braid}
A blocked braid on $n$ strings is defined to be an arrow in $\tang_n$ from $I$ to $I$ of the form
$SBR$ where $B$ is a braid in $\tang_n$.

Notice $SBR=SB'R$ does not imply that $B=B'$.
\end{dfn}

The picture of the blocked braid
$$ S(1\tn \tw)(\tw\tn 1)(1\tn \tw)(\tw\tn 1)(1\tn \tw)(\tw\tn 1)R$$
 is:
 \begin{center}
 \includegraphics[width=200mm,bb=-100 0 900 100]{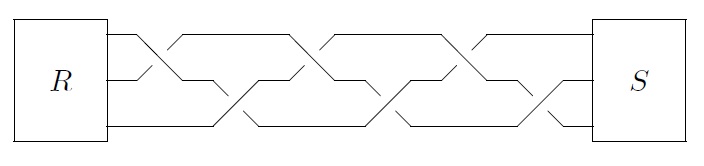}
 \end{center}

\begin{dfn}\label{def-blocked-braid-group}
The blocked-braid group on $n$ strings has as elements the blocked braids on $n$ strings. The identity of the group is $S1_{A^n}R$. The composition of $SBR$ and $SB'R$ is $SBB'R$. The inverse of $SBR$ is $SB^{-1}R$.  
\end{dfn}

We need to check that composition is well-defined. That is, to
show that if
$SBR=SB'R$ and $SCR=SC'R$ ($B,B',C,C'$ braids in $\tang_n$) then $SBCR=SB'C'R$. It clearly suffices to check this for $C=C'$.  Since the category $\tang_n$ is free with the appropriate structure, and no assumptions were made on $R$ and $S$ except their domains and codomains, then clearly $SBR=SB'R$ implies that $S'BR'=S'B'R'$ for any $S':A^n\to I$ and $R':I\to A^n$.
 Take $S'=S$ and $R'=CR$. Then $SBCR=SB'CR$. 

The group axioms are now clearly satisfied.

\begin{rem}\label{rem-composition-well-defined}
Notice the similarity of the composition operation with the sum of knots, in which a piece of each knot is removed to form the composition.
\end{rem}

\begin{dfn} The braid group $\bB_n$ on $n$ strings is generated by $n-1$ elements $\s_1, \s_2,\dots,\s_{n-1}$
and satisfies the equations 
\begin{itemize}
\item[(i)] $\s_i\s_j =\s_j\s_i$ if $i+1<j$,\\
\item[(ii)] $\s_i\s_{i+1}\s_i =\s_{i+1}\s_{i}\s_{i+1}$\\.
\end{itemize}
\end{dfn}

\begin{prop}
In $\tang_n$ (or more generally in any braided strict monoidal category) if we define 
$\s_i:A^n\to A^n$ $(i<n)$ by $\s_i=1_{A^{i-1}}\tn\tw \tn 1_{A^{n-i-1}}$ then these $\s_i$  satisfy the equations of the group $\bB_n$. This is also clearly true of the arrows $S\s_i R$ in $\BB_n$. Hence there is a surjective homomorphism from $\bB_n$ to $\BB_n$.
\end{prop}

\begin{rem}
It is convenient to use the symbol $\s_i$ in several different senses as we have done above: 1) as an element of $\bB_n$ in any of the braid groups with $n>i$, and  2) as an arrow in $\tang_n$ from $A^n$ to $A^n$ for $n>i$. In each case we will make clear by the context which of the meanings is intended for $\s_i$.
\end{rem}

 
\section{Categories of Tangled Relations}\label{sec-tangled-relations}
 
 In order to distinguish blocked-braids we describe a family of categories $\tr_G$, introduced in \cite{RSW-TangCirc-2012}. 
 
\subsection{The definition of $\tr_G$}\label{subsec-defn-trel}
We will describe a braided modification of the category $\rel$ of relations with an object with a tangle algebra structure.
\begin{definition}\label{def-trel}
Let $G$ be a group. The objects of $\tr_G$ are the formal powers of $G$, and
the arrows from $G^m$ to $G^n$
are relations $R$ from the set $G^m$ to the set $G^n$ satisfying:
\begin{itemize}
\item[1)] if $(x_1,...,x_m)R(y_1,...y_n)$ then also for all $g$ in $G$ \\
$(gx_1g^{-1}, ...,gx_mg^{-1})R(gy_1g^{-1}, ... ,gy_mg^{-1})$,
\item[2)] if $(x_1,...,x_m)R(y_1,...y_n)$ then $x_1...x_m(y_1...y_n)^{-1}\in Z(G)$ (the center of $G$).
\end{itemize}
Composition and identities are defined to be composition and identity of relations.
The tensor is defined on objects by
$G^m\tn G^n = G^{m+n}$ and on arrows by product of relations. The twist
$$\tw:G\tn G\to G\tn G$$
is the functional relation
$$(x,y)\sim (xyx^{-1},x). $$
\end{definition}

\begin{prop}\cite{RSW-TangCirc-2012}
The object $G\in \tr_G$   has a tangle algebra structure given by the arrows is the functional relation $\eta: *\sim (x,x^{-1})$, and $\eps$ is the opposite relation of $\eta$. As a consequence given any particular relations $R:G^0\to G^n$ and $S:G^n\to G^0$ in $\tr_G$ we have a representation of $\BB_n$ in $\tr_G(G^n,G^n)$.  Clearly blocked-braids which are distinct in such a representation are distinct in $\BB_n$.
\end{prop}

\subsection{The category of relations}
The usual symmetric monoidal category $\rel$ of relations  whose objects are sets, and whose arrows are relations and whose symmetry $X\times Y\to Y\times X$ is the functional relation $(x,y)\mapsto (y,x)$ has tangle algebra structures on each object  $X$: $\eta$ is the relation $\eta: *\sim (x,x)$ and $\eps$ is the opposite relation to $\eta$ \cite{CW87}. 

The representation of $\bB_n$ in $\rel(X^n,X^n)$ is that each braid $B$ goes to the corresponding 
permutation, denoted $\phi(B)$, of the factors of $X^n$. Further, given any particular relations $R:X^0\to X^n$, $S:X^n\to X^0$ (that is, subsets of $X^n$) there is a representation of the blocked-braid group which takes $SBR\in \BB_n$ to the subset of $X^0$  given by the composite of relations $S\phi(B)R:X^0\to X^n\to X^0$. It is easy to check that these representations detect the permutation of the braid in a blocked braid. Hence there is a surjective homomorphism $\BB_n$ to the symmetric group $\Sigma_n$.      
 
\section{$\BB_2$ is $\Z_2$}\label{sec-BB2}
To prove that $\BB_2=\Z_2$ it suffices to show that $S\s_1R=S\s_1^{-1}R$ in $\BB_2$. We will prove a more general result.

\begin{prop}\label{inverse-theorem}
 In $\BB_n$ the following holds:
$$S\s_{n-1}\s_{n-2}\cdots\s_1R=S\s_{n-1}^{-1}\s_{n-2}^{-1}\cdots\s_1^{-1}R.$$ Clearly also
$$S\s_{1}\s_{2}\cdots\s_{n-1}R=S\s_{1}^{-1}\s_{2}^{-1}\cdots\s_{n-1}^{-1}R.$$
\end{prop}

We give an algebraic proof but a proof which makes clear the geometric content is available in the case of $\BB_2$ in \cite{RSW-TangCirc-2012}.

\pf 
\begin{align*}
S\s_{n-1}\s_{n-2}\cdots\s_1 R &=S( 1_{A^{n-2}}\tn \tw)(1_{A^{n-3}}\tn \tw \tn 1_{A})(\cdots )(\tw\tn 1_{A^{n-2}})R\\
&=S\s_{n-1}\s_{n-2}(\cdots)(\tw\tn 1_{A^{n-2}})(1\tn \eps\tn 1_{A^{n-1}})(1_{A^2}\tn R)(\eta)\tag{duality} \\
&=S(\eps\tn 1_{A^{n}})(1_{A}\tn R\tn 1_A)(\tw\eta)\tag{naturality} \\
&=S(\eps\tn 1_{A^{n}})(1_{A}\tn R\tn 1_A)(\eta)\tag{commutativity}\\
\end{align*}
The final expression does not involve $\tw$. It is clear that repeating the above argument 
commencing with the right-hand-side reduces to the same final expression.
\qed
 
\section{$\BB_3$ has order $6$ or $12$}\label{sec-BB3}

\subsection{Equations in $\BB_3$}
Let $a=S\s_1R$ and $b=S\s_2R$ be blocked braid on three strings. 
Clearly $\BB_3$ is generated by $a$ and $b$.

The first equation satisfied by $a$, $b$ comes directly from the Yang-Baxter equation:
$$aba=S(\tw\tn 1)(1\tn \tw)(\tw\tn 1)R=S(1\tn \tw)(\tw\tn 1)(1\tn \tw)R=bab.$$

A further equation comes from the proposition of section \ref{sec-BB2}, namely
$$abba=(S\s_1\s_2R)(S\s_2\s_1R)=(S\s_1\s_2R)(S\s_2^{-1}\s_1^{-1}R)=S1R=1.$$

It is straightforward to see that the group generated by two letters $a$ and $b$ satisfying
$aba=bab$ and $abba=1$ has $12$ elements. The elements $a$ and $b$ both have order $4$ while $c=ab^{-1}$ has order $3$. Further $a^{-1}ca=c^{-1}$ hence the three element subgroup $\{1,c,c^2\}$ is normal, and in fact $$<a,b; aba=bab,abba=1>$$ is the semidirect product of $\Z_3$ by $\Z_4$ with the non-trivial action of $\Z_4$ on $\Z_3$.

Of course this does not completely identify $\BB_3$ since it could be that $\BB_3$ is the symmetric group on three elements - this would be the case if, for example also $aa=1$. 

\begin{rem} We suspect that $BB_3$ has $12$ elements but have not been able to prove it. It was pointed out in \cite{RSW-TangCirc-2012} that $\tr_G$ is unable to distinguish between the blocked braids $aa$ and $1$ in $\BB_3$.
\end{rem}

\section{Dirac's belt trick}\label{sec-dirac}

\begin{dfn}
The \emph{torsion} (through $\pi/2$) of $n$ strings, denoted $T_n$, is the braid in $\bB_n$
$$(\s_1)(\s_2\s_1)(\cdots)(\s_{n-2}\s_{n-3}\cdots\s_1)(\s_{n-1}\s_{n-2}\cdots\s_2\s_1).$$
\end{dfn}

\begin{rem} We will consider the braid 
$$(\s_1)(\s_2\s_1)(\cdots)(\s_{n-2}\s_{n-3}\cdots\s_1)(\s_{n-1}\s_{n-2}\cdots\s_2\s_1)$$
also as an element of $\bB_m$ for $m>n$ and we will still denote it $T_n$, though in
$\bB_m$ it will not be called a torsion.
\end{rem}

In \cite{RSW-TangCirc-2012} there is a geometric sketch of a proof the $ST_3^4R=S1R$ in $\BB_3$ (three strings).
The aim of this section is to prove algebraically that $(ST_nR)^4=S1R$ in $\BB_n$ ($n$ strings).
We suspect that $(ST_nR)^2$ is \emph{not} the identity in $\BB_n$ but are unable to find a proof.
In $\tr_G$ it is always the case that $ST_n^2R=S1T$. 

\begin{prop}\label{other-torsion}
The braid $T_n$ satisfies 
$$T_n= (\s_{1}\s_{2}\s_{3}\cdots \s_{n-1})(\cdots)(\s_1\s_2\s_3)(\s_1\s_2)(\s_1).$$
\end{prop}

\pf
Straightforward from the equations (ii) of braid groups.
\qed

\begin{dfn}
We denote as $Q_n$ the braid $\s_{n-1}\s_{n-2}\cdots \s_{1} $. By $Q_n^k$ we mean the $k$th power of $Q_n$. Similarly let $P_n$ be the braid $\s_{1}\s_{2}\cdots \s_{n-1} $, and let $P_n^k$ be the $k$th power of $P_n$
\end{dfn}

\begin{prop}\label{theorem4}
$(\s_{k-1}^{-1}\cdots\s_{2}^{-1}\s_1^{-1})Q_k^{k-1}=Q_{k-1}^{k-1}$.
\end{prop}

\pf
The key point in the proof is the fact that for $i<k$ 
$$\s_i^{-1}Q_{k}^{k-i}=Q_{k}^{k-i-1}Q_{k-1}.$$

The pattern will be clear from the example of $k=4$.
\begin{align*}
(\s_3^{-1}\s_2^{-1}\s_1^{-1})Q_4^3 &=\s_3^{-1}\s_2^{-1}\s_1^{-1}\s_3\s_2\s_1\s_3\s_2\s_1\s_3\s_2\s_1\\
&=\s_3^{-1}\s_2^{-1}\s_3\s_1^{-1}\s_2\s_1\s_3\s_2\s_1\s_3\s_2\s_1\\
&=\s_3^{-1}\s_2^{-1}\s_3\s_2\s_1\s_2^{-1}\s_3\s_2\s_1\s_3\s_2\s_1\\
&=\s_3^{-1}\s_2^{-1}\s_3\s_2\s_1\s_3\s_2\s_3^{-1}\s_1\s_3\s_2\s_1\\
&=\s_3^{-1}\s_2^{-1}\s_3\s_2\s_1\s_3\s_2\s_1\s_2\s_1\\
&=\s_3^{-1}\s_3\s_2\s_3^{-1}\s_1\s_3\s_2\s_1\s_2\s_1\\
&=\s_2\s_1\s_2\s_1\s_2\s_1\\
&=Q_3^3\\
\end{align*}
\qed

\begin{prop}\label{prop-Q}
The braid $T_k$ in $\bB_n$ $(n\geq k)$ satisfies $T_kT_k=Q_{k}^{k}$.
\end{prop}

\pf
The proof will be by induction on $k$. Assume the result for $k$. 
Notice that $T_{k+1}=T_k(\s_k\s_{k-1}\cdots \s_1)$.
But also by proposition \ref{other-torsion} $T_{k+1}=(\s_{1}\s_{2}\s_{3}\cdots \s_{k})T_{k}$.
Hence
$$T_{k+1}T_{k+1}= T_k(\s_k\s_{k-1}\cdots \s_1)(\s_{1}\s_{2}\s_{3}\cdots \s_{k})T_{k}.$$
But it is easy to check that $\s_i$ $(i<k)$ commutes with $(\s_k\s_{k-1}\cdots \s_1)(\s_1\s_{2}\s_{3}\cdots \s_{k})$ and hence any braid in $\bB_k$, and in particular $T_k$ commutes with
$(\s_k\s_{k-1}\cdots \s_1)(\s_{1}\s_{2}\s_{3}\cdots \s_{k})$. Hence
\begin{align*}
T_{k+1}T_{k+1}&=(\s_k\s_{k-1}\cdots \s_1)(\s_{1}\s_{2}\s_{3}\cdots \s_{k})T_kT_k\\
&=(\s_k\s_{k-1}\cdots \s_1)(\s_{1}\s_{2}\s_{3}\cdots \s_{k})Q_{k}^k\tag{inductive hypothesis}\\
&=(\s_k\s_{k-1}\cdots \s_1)(\s_{1}\s_{2}\s_{3}\cdots \s_{k})(\s_{k}^{-1}\cdots \s_{2}^{-1}\s_{1}^{-1})Q_{k+1}^k\tag{Proposition \ref{theorem4}}\\
&=(\s_k\s_{k-1}\cdots \s_1)Q_{k+1}^k\\
&=Q_{k+1}^{k+1}.\\
\end{align*}
\qed

\begin{prop}\label{prop-P}
The braid $T_k$ in $\bB_n$ $(n\geq k)$ satisfies $T_kT_k=P_{k}^{k}$.
\end{prop}

\pf
The proof follows from Proposition \ref{other-torsion} and the symmetry of the equations for the braid groups.
\qed

Finally the fact that \emph{in the blocked-braid group} four torsions equals the identity: 
\begin{thm}
In $\BB_n$, $(ST_nR)^4=S1R$.
\end{thm}

\pf
Notice that Proposition \ref{inverse-theorem} says that, in $\BB_n$, $SP_nR=SQ_n^{-1}R$. Hence in $\BB_n$
\begin{align*}
(ST_nR)^4&=(ST_n^2R)(ST_n^2R)\\
&=(SQ_n^nR)(SP_n^nR)\tag{Propositions \ref{prop-Q},\ref{prop-P}}\\
&=(SQ_n^nR)(SQ_n^{-n}R)\tag{Proposition \ref{inverse-theorem}}\\
&=SQ_n^nQ_n^{-n}R\\
&=SR.\\
\end{align*}
\qed

\section{$\BB_n$ is infinite if $n>3$}\label{sec-BB4}

\begin{thm}
$\BB_n$ is infinite if $n>3$
\end{thm}

\pf
We will show that the elements $S\s_1^kR$ $(k=0,1,2,\cdots)$ are all distinct in $\BB_n$ $(n>3)$
by looking at representations of $\BB_n$ in  $\tr_{SL_2(\Z)}$.

As the relation $R:I\to SL_2(\Z)^n$ we will take the conjugacy class of the $n$-tuple $(x,y,z,w,e,e,\cdots ,e)$
where 
$$x=
\left( \begin{array}{ccc}
1 & 1 \\
0 & 1 
\end{array} \right),
y=
\left( \begin{array}{ccc}
1 & 0 \\
1 & 1 
\end{array} \right),
z=
\left( \begin{array}{ccc}
1 & 0 \\
-1 & 1 
\end{array} \right),
w=
\left( \begin{array}{ccc}
1 & -1 \\
0 & 1 
\end{array} \right),
e=
\left( \begin{array}{ccc}
1 & 0 \\
0& 1 
\end{array} \right).
$$ 
Notice that $xyzwe^{n-4}=e\in Z(SL_2(\Z))$ and hence the relation is in $\tr_{SL_2(\Z)}$.
We will now consider just the case $n=4$ because the argument extends trivially to the cases $n>4$.
We will show that the conjugacy classes $S_k=\s_1^kR$ $(k=0,1,2,\cdots)$ are pairwise disjoint and hence
$S_k\s_1^lR=\emptyset$ except when $k=l$ in which case it is the one point set.

The conjugacy classes $S_k$ $(k=0,1,2,\dots)$ are the conjugacy classes of the $4$-tuples
\begin{align*}
(x,y,z,w), (xyx^{-1},x,z,w), &(xyxy^{-1}x^{-1},
xyx^{-1},z,w),\\
&(xyxyx^{-1}y^{-1}x^{-1},xyxy^{-1}x^{-1},z,w),\cdots\\
\end{align*}     
If $(u_1,v_1,z,w)$, $(u_2,v_2,z,w)$ are two different elements in this list which generate the same conjugacy classes there must exist an element $g\in SL_2(\Z)$ such that $g^{-1}u_1g=v_1$, 
$g^{-1}u_2g=v_2$, $g^{-1}zg=z$, $g^{-1}wg=w$. But $g^{-1}zg=z$, $g^{-1}wg=w$ imply that $g$ is either the identity matrix $e$ or $-e$.  Then $g^{-1}u_1g=v_1$, $g^{-1}u_2g=v_2$ imply that $u_1=v_1$ and $u_2=v_2$.
However a straightforward calculation of the sequence of matrices
 $$x,xyx^{-1},xyxy^{-1}x^{-1},xyxyx^{-1}y^{-1}x^{-1},xyxy^{-1}x^{-1},\cdots$$ 
shows that no repetition occurs - the maximum absolute value of the entries of the $i$th matrix in the sequence is strictly increasing as $i$ increases.
\qed  
 

\end{document}